\documentclass[12pt]{amsart}

\newtheorem{theorem}{Theorem}[section]
\newtheorem{lemma}[theorem]{Lemma}
\newtheorem{corollary}[theorem]{Corollary}
\theoremstyle{definition}   
\newtheorem{definition}[theorem]{Definition}

\theoremstyle{remark}
\newtheorem{remark}[theorem]{Remark}

\newtheorem{proposition}[theorem]{Proposition}
\numberwithin{equation}{section}

\usepackage{times}
\usepackage{enumerate}
\usepackage{tfrupee}
\usepackage{tikz}
\usetikzlibrary{chains,fit}
\usetikzlibrary{shapes,snakes}
\usetikzlibrary{graphs}
\usepackage{graphicx,adjustbox}
\usepackage{hyperref}
\usepackage{amsmath,amssymb}

\newcommand{\lcm}{{\rm \mbox{lcm}}}
\newcommand{\LT}{{\rm \mbox{Lt}}}

\usepackage{amscd}
\usepackage{graphicx}
\usepackage[all]{xy}

\title[Primary decomposition and normality of certain determinantal ideals]
{Primary decomposition and normality of certain determinantal ideals}
\author{
Joydip Saha
\and
Indranath Sengupta
\and
Gaurab Tripathi
}
\date{}

\address{\small \rm  Discipline of Mathematics, IIT Gandhinagar, Palaj, Gandhinagar, 
Gujarat 382355, INDIA.} 
\email{saha.joydip56@gmail.com}
\thanks{The first author thanks SERB for the Post-doctoral fellowship under the 
research grant EMR/2015/000776.}

\address{\small \rm  Discipline of Mathematics, IIT Gandhinagar, Palaj, Gandhinagar, 
Gujarat 382355, INDIA.}
\email{indranathsg@iitgn.ac.in}
\thanks{The second author is the corresponding author, who is supported by the 
research project EMR/2015/000776 sponsored by the SERB, Government of India.}

\address{\small \rm Department of Mathematics, Jadavpur University, Kolkata,
WB 700 032, India.} 
\email{gelatinx@gmail.com}
\thanks{The third author thanks CSIR for the Senior Research Fellowship.}

\date{}



\begin{document}
\begin{abstract}
In this paper we study primality and primary decomposition of certain ideals which are 
generated by homogeneous degree $2$ polynomials and occur naturally from determinantal 
conditions. Normality is derived from these results.
\end{abstract}

\maketitle

\section{Introduction}
Let $K$ be a field. Let $\{x_{ij}; \, 1\leq i \leq m, \, 1\leq j \leq n\}$, 
$\{y_{j}; \, 1\leq j \leq n\}$ be indeterminates over $K$, so that 
$R=K[x_{ij}, y_{j}]$ denotes the polynomial algebra over $K$. Let $X$ 
denote an $m\times n$ matrix such that its entries belong to the ideal 
$\langle \{x_{ij}; \, 1\leq i \leq m, \, 1\leq j \leq n\}\rangle$. Let 
$Y=(y_{j})_{n\times 1}$ be the generic $n\times 1$ column matrix. Let 
$I_{1}(XY) = \langle g_{1}, \ldots , g_{m}\rangle$ denote 
the ideal generated by the $1\times 1$ minors or the entries of the 
$m\times 1$ matrix $XY$. 
\medskip

Ideals of the form $I_{1}(XY)$ have been studied by several authors, see 
\cite{hocheagon}, \cite{herzog}, \cite{hunu1}, \cite{hunu2}. Subsequently, 
ideals of the form $I_{1}(XY)+J$, where $J$ is also determinantal appeared 
in the paper \cite{johnson}. Our aim is to construct explicit primary 
decompositions of ideals of the form $I_{1}(XY)$ by constructive techniques. 
The cases we consider are the following:
\begin{enumerate}
\item $m=n$ and $X$ is generic or generic symmetric;
\item $m=n+1$ and $X$ is generic. 
\end{enumerate}
The articles \cite{hocheagon}, \cite{hunu1}, \cite{hunu2} have also discussed 
certain issues like primality of these ideals, when $X$ is generic. However, our 
techniques only use the information about Gr\"{o}bner basis for these 
ideals from \cite{sst} and the notion of complete irreducibility from 
\cite{fer}. Moreover, this study has been used to prove that $I_{1}(XY)$ 
is normally torsionfree in \ref{torfree} and hence normal in \ref{normality}. 
The case when $X$ is generic skew-symmetric turns out to be the most 
challenging one. The best result we could obtain is in \ref{primdec2}. 
The primary decompositions of $I_{1}(XY)$ is not known when $X$ is an 
$n\times n$ generic skew-symmetric matrix.
\medskip

First we note that the ideal $I_{1}(XY)$ is 
not a prime ideal if $m=n$ and $X$ is one of the above. Let 
us prove this for $m = n$ and $X$ generic. A similar proof follows 
if $X$ is symmetric. Let $\Delta=\det(X)$. It is easy to see that 
\begin{eqnarray*}
\Delta \cdot y_{n} & = & (\sum_{j=1}^{n}A_{jn}x_{jn})y_{n}, \, 
A_{jn} \, {\rm being\, the} \, (j,n)- \, {\rm cofactor\, of} 
\, X;\\
{} &  = & \sum_{j=1}^{n}A_{jn}\left(\sum_{k=1}^{n}x_{jk}y_{k}\right) 
 - \sum_{j=1}^{n}A_{jn}\left(\sum_{k\neq n}x_{jk}y_{k}\right)\\
{} & = & \sum_{j=1}^{n}A_{jn}g_{j},
\end{eqnarray*} 
since $\sum_{j=1}^{n}A_{jn}\left(\sum_{k\neq n}x_{jk}y_{k}\right) 
 = \sum_{k\neq n}\left(\sum_{j=1}^{n}A_{jn}x_{jk}\right)y_{k} = 0$. 
Therefore  $\Delta\cdot y_{n}\in I_{1}(XY)$, but $\Delta \notin I_{1}(XY)$ 
and $y_{n}\notin I_{1}(XY)$. Similar argument as above shows that if 
$m=n$ and $X$ is generic symmetric, then $I_{1}(XY)$ is not a prime ideal. 
Primary decompositions of $I_{1}(XY)$ for the cases $m=n$ and $m=n+1$ are given in 
section 5.
\medskip

In the case when $m=n$ and $X$ is generic skew-symmetric, it is easy to see that 
the ideal $\langle g_{1}, \ldots , g_{n-1} \rangle $ is not a prime ideal 
and the sequence $g_{1}, \ldots , g_{n}$ is not regular, for  
$y_{n}g_{n} = (-y_{1})g_{1}+ (-y_{2})g_{2}+\ldots+ (-y_{n-1})g_{(n-1)}$, 
but $y_{n},g_{n}\notin \langle g_{1}, \ldots , g_{n-1} \rangle$.
\medskip

The main theorems proved in this paper are the following: 
\medskip

\begin{theorem}\label{primdec1}
Let $X$ denote an $m\times n$ matrix such that its entries belong to the ideal 
$\langle \{x_{ij}; \, 1\leq i \leq m, \, 1\leq j \leq n\}\rangle$ and $Y=(y_{j})_{n\times 1}$ denote the generic $n\times 1$ column matrix. Let 
$\mathcal{I}_{t} = \langle g_{1}, \ldots , g_{t}\rangle$, with $1\leq t\leq m$. 
\begin{enumerate}
\item Let $m=n$ and $X = (x_{ij})$ be generic or generic symmetric. 
\begin{enumerate}
\item[(i)] $g_{1}, \ldots , g_{n}$ is a regular sequence in $R$.
\item[(ii)] For every $1\leq t\leq n-1$, the ideal $\mathcal{I}_{t}$ is a prime ideal 
in $R$.
\item[(iii)] The primary decomposition of the ideal $\mathcal{I}_{n} = I_{1}(XY)$ is  
given by 
$$\mathcal{I}_{n}=\langle y_{1},\cdots,y_{n}\rangle\cap \langle g_{1},\cdots, g_{n},\Delta\rangle,$$
where $\Delta$ denotes the determinant of $X$.
\end{enumerate}

\item Let $m=n+1$ and $X = (x_{ij})$ be generic. The primary decomposition of the ideal 
$\mathcal{I}_{n+1} = I_{1}(XY)$ is given by 
$$\mathcal{I}_{n+1}= \langle y_{1},\cdots,y_{n} \rangle \cap \langle g_{1},\cdots, g_{n},\Delta_{1},\cdots,\Delta_{n+1}\rangle,$$ 
where $\Delta_{i}$  denotes the determinant of the $n\times n$ matrix formed by removing the 
$i$-th row of the matrix $X$. 
\end{enumerate}
\end{theorem}
\medskip

\begin{remark}\label{corollary}
It follows from part (1), statement (ii) of the above theorem that if $m< n$ and 
$X = (x_{ij})$ is generic, then $I_{1}(XY)$ is a prime ideal.
\end{remark}
\medskip

\begin{theorem}\label{primdec2}
Let ${\rm ch}(K)\neq 2$. Let $m=n$ and $X = (x_{ij})$ be 
generic skew-symmetric. Let 
$\mathcal{J}_{t} = \langle g_{1}, \ldots , g_{t}\rangle$, with $1\leq t\leq n$. 
\begin{enumerate}
\item[(i)] $g_{1}, \ldots , g_{n-1}$ is a regular 
sequence in $R$. 
\item[(ii)] For every $1\leq t\leq n-2$, the ideal 
$\mathcal{J}_{t}=\langle g_{1},\ldots,g_{t}\rangle$ is a prime ideal in $R$.
\end{enumerate}
\end{theorem}
\medskip

It is not difficult to prove that 
$g_{1}, \ldots , g_{n}$ (respectively  $g_{1}, \ldots , g_{n-1}$) 
is a regular sequence, if we choose a monomial order on $R$ suitably. 
For proving primality and primary decomposition we need techniques developed 
in \cite{fer}; see sections 3 and 4. Primality can also be proved by 
geometric arguments, which is perhaps more natural for ideals of this 
form. However, we could not figure out better technique for describing 
primary decomposition for these ideals.
\bigskip

\section{Regular sequence}
We begin with some rudiments of {\it Gr\"obner bases}, which is used 
as the main tool for proving most of the statements. Let 
$K[x_{1},\ldots,x_{n}]$ be the polynomial ring over $K$, 
with a monomial order $>$. For 
$0 \neq f\in K[x_1,\ldots,x_n]$, let ${\rm Lm}_{>}(f)$, ${\rm Lt}_{>}(f)$, 
${\rm Lc}_{>}(f)$ denote the {\it leading monomial}, {\it leading term} and 
{\it leading constant} of $f$ respectively. 
We simply write ${\rm Lt}(f)$, ${\rm Lc}(f)$, ${\rm Lm}(f)$, 
when no confusion is likely to occur. For $0\neq f, g\in K[x_{1},\ldots,x_{n}]$, 
the {\it S - polynomial} of $f$ and $g$, denoted by $S(f, g)$, 
is the polynomial 
$$S(f, g) := \frac{\lcm({\rm Lm}(f), {\rm Lm}(g))}{{\rm Lt}(f)}\cdot f - 
\frac{\lcm({\rm Lm}(f), {\rm Lm}(g))}{{\rm Lt}(g)}\cdot g\,.$$ 
Given $G = \{g_{1},\ldots, g_{t}\}\subseteq K[x_1,\ldots,x_n]$ and 
$f\in K[x_{1},\ldots,x_{n}]$, we say that 
$f$ {\it reduces to zero modulo} $G$, 
denoted by $f\to_{G} 0$, if $f$ can be written as 
$f = \sum_{i=1}^{t} a_{i}g_{i}\,$, 
such that $\,{\rm Lm}(f) \geq {\rm Lm}(a_{i}g_{i})$, whenever 
$\,a_{i}g_{i}\neq 0$. Buchberger's Criterion says that, for an ideal $I$ 
in $K[x_{1},\ldots,x_{n}]$ and a generating set $G = \{g_{1},\ldots, g_{t}\}$ 
for $I$, $G$ is a Gr\"obner basis for $I$ iff 
$S(g_{i}, g_{j}) \to_{G} 0$, for every $i\neq j$. 
The following simple lemma is useful in almost any computation 
involving Gr\"{o}bner bases. 
\medskip

\begin{lemma} 
Let $G = \{g_{1},\ldots, g_{t}\}\subseteq K[x_{1},\ldots,x_{n}]$ and let 
$f, g\in G$ be non-zero with ${\rm Lc}(f) = {\rm Lc}(g) = 1$, and 
$\gcd(\,{\rm Lm}(f), \,{\rm Lm}(g)\,) = 1$. Then, 
\begin{enumerate}
\item $S(f, g) \,=\,  {\rm Lm}(g).\,f \,-\, {\rm Lm}(f).\,g$\,.

\item $S(f, g) \,=\,  -(g - {\rm Lm}(g)).\,f \,+\, (f - {\rm Lm}(f)).\,g\, 
\longrightarrow_{G} 0$. 

\end{enumerate}
\end{lemma}
\medskip

In order to show that a set $\{g_{1},\ldots, g_{t}\}\subseteq 
K[x_{1},\ldots,x_{n}]$ is a Gr\"{o}bner 
basis, one always intends to choose the monomial order in such a way 
that most of the pairs of the leading monomials ${\rm Lm}(g_{i})$ 
are mutually coprime. This would be visible in \ref{disjoint} and also 
in \ref{regseq}, \ref{grob}. However, the cases for generic and generic skew-symmetric are quite different. The 
case when $X$ is generic is easier to handle and a 
diagonal term order is often useful, which is probably due to the fact that 
there is no relation among the entries of $X$. The case of 
generic symmetric is also not so different from the generic case. On 
the other hand, when $X$ is generic skew-symmetric, the presence of $0$ and other 
relations in the entries of the matrix $X$ affect the entries of 
$I_{1}(XY)$ and make the situation far more complicated.
\medskip
 
\begin{lemma}\label{disjoint}
Let $h_{1},h_{2}\cdots, h_{n}\in R$ be such that with respect to a suitable 
monomial order on $R$, the leading terms of them are mutually coprime. 
Then, $h_{1},h_{2}\cdots, h_{n}$ is a regular sequence in $R$.
\end{lemma}

\proof. The element $h_{1}$ is a regular element in $R$, since $R$ is a domain and 
$h_{1}\neq 0$. By induction we assume that for $k\leq n-1$, 
$\{h_{1},h_{2}\cdots, h_{k}\}$ forms a regular sequence in $R$. We note that 
the set $\{h_{1},h_{2}\cdots, h_{k}\}$ is a Gr\"{o}bner basis for the ideal 
$J=\langle h_{1},h_{2}\cdots, h_{k}\rangle$, since $\gcd({\rm \LT}(h_{i}), {\rm \LT}(h_{j}))=1$ for every $i\neq j$. Let 
$gh_{k+1}\in J$. Then 
$\LT(g)\LT(h_{k+1})$ must 
be divisible by $\LT(h_{i})$ for some $1\leq i\leq k$. But, 
$\gcd(\LT(h_{i}), \LT(h_{k+1}))=1$, and hence $\LT(h_{i})$ divides 
$\LT(g)$. Let $r = g - \dfrac{\LT(g)}{\LT(h_{i})}h_{i}$. If $r=0$, 
then $g\in J$. If $r\neq 0$, then $\LT(r)<\LT(g)$ and $rh_{k+1}\in J$. We follow the same argument with $rh_{k+1}$.\qed
\medskip

\begin{theorem}\label{regseq}
\begin{enumerate}
\item[(i)] Let $X = (x_{ij})_{m\times n}$ be either generic with $m\leq n$ or generic symmetric with $m=n$. Then $g_{1}, \ldots , g_{m}$ is a regular sequence in $R$.
\medskip

\item[(ii)] Let ${\rm ch}(K)\neq 2$. Let $m=n$ and $X = (x_{ij})$ be generic skew-symmetric, then $g_{1}, \ldots , g_{n-1}$ is a regular 
sequence in $R$. 
\end{enumerate}
\end{theorem}

\proof To prove (i), we choose the lexicographic monomial order on $R$ given by
$$x_{11}> x_{22}> \cdots >x_{mm}; \quad x_{ij}, y_{j} < x_{mm},$$
for every $1 \leq i \leq m$, $1\leq j \leq n$ and $i\neq j$.
\medskip

\noindent To prove (ii), we choose the lexicographic monomial order on $R$ given by
$$x_{12}< x_{23}< \cdots <x_{(n-1)n}; \quad x_{ij}, y_{j} < x_{12},$$ 
for every $1\leq i \leq n-1$, $2\leq j \leq n$, $i<j$ and $i\neq j-1$.
\medskip

\noindent The leading terms of the polynomials $g_{1}, \ldots , g_{m}$ 
in the case (i) and $g_{1}, \ldots , g_{n-1}$ in the case (ii) are 
mutually coprime with respect to the respective monomial orders defined 
above. We now apply Lemma \ref{disjoint} to prove the statement.\qed
\bigskip

\section{Primality of $\mathcal{I}_{t}$ and $\mathcal{J}_{t}$}
We prove primality of the ideals $\mathcal{I}_{t}$ and $\mathcal{J}_{t}$ 
defined in Theorems \ref{primdec1} and \ref{primdec2} using Theorem 2.5 in \cite{fer}. 
Let us first recall the notion of \textit{complete irreducibility} from \cite{fer}. 
\medskip

\noindent\textbf{Complete Irreducibility.} 
Let $P$ be a commutative ring with identity. Let $\mathfrak{a}$ be a prime ideal of $P$. Let 
$$\Gamma_{\mathfrak{a}}:= \{f\in P[x]\mid arf-fra \in P[x]; \forall r \in P, \delta f\neq 0 ,\ a= \text{lc}(f)\notin \mathfrak{a}\},$$ 
where $\delta f$ denotes the degree of $f$ and $\text{lc}$ denotes the leading coefficient of $f$, with respect to the indeterminate $x$. Given $f\in\Gamma_{\mathfrak{a}}$, let
$$[\mathfrak{a},f]:=\{g\in P[x]\mid g\langle a^{e}\rangle\subset 
\mathfrak{a}[x]+\langle f\rangle\, \, \text{for\, some\, integer}\, e\geq 0\}.$$ 
A polynomial $f\in\Gamma_{\mathfrak{a}}$ is $\Gamma_{\mathfrak{a}}$ 
\textit{completely irreducible} if the following criteria holds: \, If 
there exist $b\in P$, $g\in \Gamma_{\mathfrak{a}}$, $h\in P[x]$, such that 
$fb\notin P[x]$ and $fb-gh\in P[x]$ then $\delta g=\delta f$.
\medskip

Let $Q=P[x_{1},\ldots,x_{n}]$. For $i=1,\ldots, n$, let 
$f_{i}\in P[x_{1},\ldots,x_{i-1}][x_{i}]$, with 
$a_{i}=\text{lc}(f_{i})\in P[x_{1},x_{2},\cdots,x_{i-1}]$, with respect to 
the indeterminate $x_{i}$. Let
\begin{eqnarray*}
[\mathfrak{a},f_{1},\cdots,f_{n}] & = & 
\{g\in Q \mid g\langle a_{1}\rangle^{e_{1}}\cdots \langle a_{n}\rangle ^{e_{n}}\subset 
\mathfrak{a}[x_{1},\cdots,x_{n}]+ \langle f_{1},\cdots,f_{n}\rangle,\\                                                                                                    
{} & {} & \quad\quad\quad\quad\quad\quad \text{for\, nonnegative\, integers} 
\, e_{1}, \ldots , e_{n}\}.
\end{eqnarray*} 
If $\mathfrak{a} = \langle 0\rangle $ is prime then $[\langle 0\rangle,f_{1},\cdots,f_{n}]$ is 
written as $[f_{1},\cdots,f_{n}]$. Therefore, 
$$[f_{1},\cdots,f_{n}]=
\{g\in Q \mid g\langle a_{1}\rangle^{e_{1}}\cdots \langle a_{n}\rangle^{e_{n}}\subset 
\langle f_{1},\cdots,f_{n}\rangle,\, \text{for\, integers\,} \, e_{i}\geq 0\}.$$ 
The sequence $(f_{1},f_{2},\cdots,f_{n})$ defined above is said to be 
\textit{completely irreducible (mod $\mathfrak{a}$)} if $f_{1}$ is $\Gamma_{\mathfrak{a}}$ 
completely irreducible and $f_{i+1}$ is $\Gamma_{\mathfrak{a}_{i}}$ completely irreducible as a 
polynomial in $x_{i+1}$, where $\mathfrak{a}_{0} = \mathfrak{a}$ and 
$\mathfrak{a}_{i}= [\mathfrak{a}_{i-1},f_{i}]$, for every $0\leq i\leq n-1$. 
We now state Theorem 2.5 in \cite{fer} which is the main tool for 
proving primality of $\mathcal{I}_{t}$ and $\mathcal{J}_{t}$:
\medskip

\noindent\textbf{Theorem} (2.5; \cite{fer}). Suppose that $(f_{1},\cdots,f_{n})\in\mathcal{F}$ is a completely irreducible system (mod $\mathfrak{a}$). 
The ideal $\mathfrak{b}=[\mathfrak{a},f_{1},\cdots,f_{n}]$ is a prime ideal 
of $Q$, such that $\mathfrak{b}\cap P = \mathfrak{a}$.
\medskip
 
\begin{lemma}\label{compirr}
\begin{enumerate}
\item[(i)] Let $m=n$ and $X = (x_{ij})$ be generic or generic symmetric. 
The sequence $(g_{1}, \ldots , g_{n})$ is completely irreducible 
(mod $\langle 0\rangle$) and the ideal $[g_{1},\cdots,g_{n}]$ is a prime ideal.
\medskip

\item[(ii)] Let ${\rm ch}(K)\neq 2$. Let $m=n$ and $X = (x_{ij})$ be generic skew-symmetric. The sequence 
$(g_{1}, \ldots , g_{n-1})$ is completely irreducible 
(mod $\langle 0\rangle$) and the ideal $[g_{1},\cdots,g_{n-1}]$ is a prime ideal.
\end{enumerate}
\end{lemma}

\proof We prove the statement only for the first case, that is if $m=n$ and 
$X = (x_{ij})_{m\times n}$ is generic. The proofs for the other two cases 
are similar. 
\medskip

Let $P:=K[x_{ij},y_{i}\mid 1\leq i\leq m, 1\leq j\leq n-1]$ and $\mathfrak{a}_{0}:=(0)$. Then $g_{1}\in P[x_{1n}]=:P_{1}$, $g_{2}\in P_{1}[x_{2n}]=:P_{2}$, and so on 
$g_{m}\in P_{m-1}[x_{mn}] = R$. We show that the sequence $(g_{1},\ldots,g_{m})$ is completely irreducible (mod $\langle 0\rangle$). We have 
$\Gamma_{\mathfrak{a}_{0}}=\Gamma_{\langle 0\rangle}=\{f\in P_{1}\mid \delta f\neq 0\, \textrm{and}\, \text{lc}(f)\neq 0\}$. It is clear that 
$g_{1}\in \Gamma_{\langle 0\rangle}$ and we show that $g_{1}$ 
is $\Gamma_{\langle 0\rangle}$ irreducible. Suppose that 
$b\in P$, $g\in \Gamma_{0}$, $h\in P_{1}$, with $b\neq 0$ and $g_{1}b-hg=0$. 
Now $\delta g\geq 1$ as a polynomial of $x_{1n}$, since $g\in \Gamma_{0}$. If $\delta (g)>1$ then the degree of $hg$ as polynomial in $x_{1n}$ is greater than one, on the other hand the degree of $g_{1}b$ as polynomial in $x_{1n}$ is exactly one; which is a contradiction. Therefore, $\delta (g)$ must be equal to $1$.
\medskip

By induction let us assume that the sequence $(g_{1},\ldots ,g_{i-1})$ is a completely irreducible system (mod($0$)). Then $\mathfrak{a}_{i-1}:= [\langle 0\rangle,g_{1},\ldots, g_{i-1}]$ is a prime ideal by Theorem 2.5 in \cite{fer}. We first show that $y_{n}\notin \mathfrak{a}_{i-1}$, for all $i\geq 1$. If $i=1$, 
$\mathfrak{a}_{0}=\langle 0\rangle$; hence $y_{n}\notin \mathfrak{a}_{0}$. 
Let us assume that it holds for $i=t-1$. We know that 
$$\mathfrak{a}_{t}=[\mathfrak{a}_{t-1},g_{t}]=\{g\in P_{t}\mid \exists\, e\geq 0\,\textrm{with}\,\,g\langle y_{n}^{e}\rangle\subset P_{t}g_{t}+\mathfrak{a}_{t-1}[x_{tn}]\}.$$ 
If $y_{n}\in \mathfrak{a}_{t}$, then 
$y_{n}^{e}\in P_{t}g_{t}+\mathfrak{a}_{t-1}[x_{tn}]$, for some $e\geq 1$. We can write  
$y_{n}^{e} = p.g_{t} + q$, for some $p\in P_{t}$ and $q\in \mathfrak{a}_{t-1}[x_{tn}]$. On 
substituting $x_{ij} = 0$ in the above expression we get $y_{n}^{e} = c$ and 
$c\in \mathfrak{a}_{t-1}$; which contradicts the induction hypothesis. Therefore $y_{n}\notin \mathfrak{a}_{t}$. Given that $\delta g_{i}=1$ and $\text{lc}(g_{i})=y_{n}\notin \mathfrak{a}_{i-1}$ 
as a polynomial of $x_{in}$, we have $g_{i}\in \Gamma_{\mathfrak{a}_{i-1}}$. We now 
show that $g_{i}$ is $\Gamma_{\mathfrak{a}_{i-1}}$ irreducible. Suppose that $b\in P_{i-1}$, 
$g\in \Gamma_{\mathfrak{a}_{i-1}}$, $h\in P_{i-1}[x_{in}]$, such that 
$g_{i}b\notin \mathfrak{a}_{i-1}[x_{in}]$ and $g_{i}b-hg\in  \mathfrak{a}_{i-1}[x_{in}]$. Let $g_{i}b-hg=\sum_{p=0}^{t} c_{p}x_{in}^{p}$, where $c_{p}\in \mathfrak{a}_{i-1}$. Let us write $g_{i}=x_{in}y_{n}+c$, where $c\in P_{i-1}$, $h=\sum_{p=0}^{l} b_{p}x_{in}^{p}$ and $g=\sum_{p=0}^{r} a_{p}x_{in}^{p}$. Since $g\in \Gamma_{\mathfrak{a}_{i-1}}$ 
we have $a_{r}\notin \mathfrak{a}_{i-1}$. Since $g_{i}b\notin \mathfrak{a}_{i-1}[x_{in}]$, there exist $b_{p}$ such that $b_{p}\notin \mathfrak{a}_{i-1}$. Without loss of 
generality we may assume that $b_{l}\notin \mathfrak{a}_{i-1}$, otherwise we take 
$g_{i}b-(h-\sum_{p=k+1}^{l}b_{p}x_{in}^{p})g\in \mathfrak{a}_{i-1}[x_{in}]$, 
where $k: = \text{max}\{p\mid b_{p}\notin P_{i-1}\}$. Consider the equation,
$$b(x_{in}y_{n}+c)=(\sum_{p=0}^{l} b_{p}x_{in}^{p})(\sum_{p=0}^{r} a_{p}x_{in}^{p})+(\sum_{p=0}^{t} c_{p}x_{in}^{p}).$$
Now $\mathfrak{a}_{i-1}$ is a prime ideal, $a_{r}\notin \mathfrak{a}_{i-1}$ and 
$b_{l}\notin \mathfrak{a}_{i-1}$ imply that  $b_{l}a_{r}\notin \mathfrak{a}_{i-1}$; 
while each coefficient of $\sum_{p=0}^{t} c_{p}x_{in}^{p}$ is in $\mathfrak{a}_{i-1}$. Therefore, no term of $\sum_{p=0}^{t} c_{p}x_{in}^{p}$ can cancel 
with $a_{r}b_{l}x_{in}^{r+s}$. Equating degree as a polynomial of $x_{in}$, we have $r+s=1$. Therefore, we must have $r=1$, since $r\geq 1$ and $s\geq 0$. Hence, $g_{1},\ldots ,g_{i}$ is completely irreducible system mod$\langle 0\rangle $. Primality of the ideal follows 
from Theorem 2.5 in \cite{fer}. \qed
\medskip

\begin{corollary}
\begin{enumerate}
\item[(i)] Let $m=n$ and $X = (x_{ij})$ be generic or generic symmetric. 
The sequence $(g_{1}, \ldots , g_{t})$ is completely irreducible 
(mod $\langle 0\rangle$) and the ideal $[g_{1},\cdots,g_{t}]$ is a prime ideal, 
for every $t=1, \ldots , n$.
\medskip

\item[(ii)] Let ${\rm ch}(K)\neq 2$. Let $m=n$ and $X = (x_{ij})$ be generic skew-symmetric. The sequence 
$(g_{1}, \ldots , g_{t})$ is completely irreducible 
(mod $\langle 0\rangle$) and the ideal $[g_{1},\cdots,g_{t}]$ is a prime ideal, 
for every $t=1, \ldots , n-1$.
\end{enumerate}
\end{corollary}

\proof The proof follows from the definition of complete irreducibility and 
Theorem 2.5 in \cite{fer}.\qed
\medskip

\begin{theorem}\label{idealbox}
\begin{enumerate}
\item[(i)] Let $m=n$ and $X = (x_{ij})$ be generic or generic symmetric. Then, 
$[g_{1},\cdots,g_{t}] = \langle g_{1},\cdots,g_{t}\rangle$, \, for \, every \, 
$t=1, \ldots , n-1$.
\medskip

\item[(ii)] Let ${\rm ch}(K)\neq 2$. Let $m=n$ and $X = (x_{ij})$ be generic skew-symmetric. Then, 
$[g_{1},\cdots,g_{t}] = \langle g_{1},\cdots,g_{t}\rangle$, \, for \, every \, 
$t=1, \ldots , n-2$.
\end{enumerate}
\end{theorem}

\proof The proofs for the generic and generic symmetric cases would require 
the monomial order defined in the proof of Theorem \ref{regseq}; part (i). 
The proof for the generic skew-symmetric case is similar, with the only exception that 
it would require the monomial order defined in the proof of Theorem \ref{regseq}; 
part (ii). We prove the statement only for the first case, that is if $m=n$ and 
$X = (x_{ij})$ is generic. The proofs for the other two cases are similar. 
\medskip

Let $1\leq t\leq n-1$. It is clear from the definition that 
$\langle g_{1},\cdots,g_{t}\rangle \subseteq [g_{1},\cdots,g_{t}]$. Let 
$g\in [g_{1},\cdots,g_{t}]$. Then, $g\cdot y_{n}^{e}\in \langle g_{1},\cdots,g_{t}\rangle$. 
We know that $\LT(g_{i}) = x_{ii}y_{i}$, for every $1\leq i\leq t < n$, 
with respect to the monomial order chosen in the Theorem \ref{regseq}. 
We also know that the leading terms of $g_{1},\cdots,g_{t}$ are mutually coprime and 
therefore they form a Gr\"{o}bner basis for the ideal $\mathcal{I}$ with respect 
to the said monomial order on $R$. It is clear that 
$y_{n}$ does not divide $\LT(g_{i})$ for every $1\leq i\leq t$ and hence 
$\LT(g_{i})\mid \LT(g)$ for some $1\leq i\leq t$. We can write 
$g = q_{1}g_{1} + \cdots + q_{t}g_{t} + r$, where $q_{1}, \ldots , q_{t}\in R$ and 
$r$ is the remainder. Therefore, 
$$r\cdot y_{n}^{e} = g\cdot y_{n}^{e} - (q_{1}g_{1} + \cdots + q_{t}g_{t})\cdot y_{n}^{e}\in \mathcal{I}.$$
If $r\neq 0$ then $\LT(g_{i})\mid \LT(r\cdot y_{n}^{e})$ for some $1\leq i\leq t$. 
Therefore, $\LT(g_{i})\mid \LT(r)$ for some $1\leq i\leq t$, which contradicts the 
fact that $r$ is the remainder. Hence, $r=0$ and this proves that $g\in \mathcal{I}$.\qed
\medskip

\begin{remark}
In the case when $m = n$ and $X$ is generic or generic symmetric, we have mentioned 
in the introduction that $\Delta\cdot y_{n}\in \langle g_{1},\cdots,g_{n}\rangle$. Therefore, 
$\Delta\in [g_{1},\cdots,g_{n}]$ but $\Delta \notin \langle g_{1},\cdots,g_{n}\rangle$, proving that 
$[g_{1},\cdots,g_{m}] \neq \langle g_{1},\cdots,g_{m}\rangle$. 
Similarly, in the case $m = n$ and $X$ is generic skew-symmetric, we have seen in the 
introduction that $y_{n}g_{n} = (-y_{1})g_{1}+ (-y_{2})g_{2}+\ldots+ (-y_{n-1})g_{(n-1)}$. Therefore, $g_{n}\in [g_{1},\cdots,g_{n-1}]$ but 
$g_{n} \notin \langle g_{1},\cdots,g_{n-1}\rangle$, proving that 
$[g_{1},\cdots,g_{n-1}] \neq \langle g_{1},\cdots,g_{n-1}\rangle$. 
\end{remark}
\bigskip

\section{Primary decomposition of $I_{1}(XY)$}

\begin{theorem}\label{grob}
Let $I= \langle g_{1},\dots, g_{n}, \Delta\rangle$ and $\mathfrak{G}= \left(\mathcal{G}\setminus G_{n}\right)\cup\{\Delta\}$, where 
$\mathcal{G}$ is a Gr\"obner basis for the ideal 
$\langle g_{1},\cdots, g_{n}\rangle$ described in \cite{sst}. 
Then $\mathfrak{G}$ is a Gr\"obner basis for $I$, with respect to the lexicographic monomial order given by $y_{1}>\cdots >y_{n}>x_{11}>x_{12}>\cdots x_{n,(n-1)}>x_{n,n}$ on $R$. 
\end{theorem}

\proof
We note that $G_{n} = \{\Delta y_{n}\}$ and $\LT(\Delta y_{n}) = \LT(\Delta)y_{n}$. 
Hence $\LT(\mathfrak{G})=\LT(\mathcal{G}\setminus \{\Delta y_{n}\})\cup\{\LT(\Delta)\}$. We apply Buchberger's criterion. Let $f,g\in\mathfrak{G}$. Either $f$ or $g$ must 
belong to $\mathcal{G}$ since $\mathfrak{G}$ differs from $\mathcal{G}$ only by 
a single element. We consider two cases separately.
\medskip

Suppose that $f, g\in \mathcal{G}$. Then, the S-polynomial 
$S(f,g)\longrightarrow_{\mathcal{G}} 0$, 
since $\mathcal{G}$ forms a Gr\"obner basis for the ideal 
$\langle g_{1},\cdots, g_{n}\rangle$.  Hence $S(f,g)\longrightarrow_{\mathfrak{G}} 0$, since $\mathcal{G}$ and $\mathfrak{G}$ differ by a single element and 
$\LT(\Delta y_{n}) = \LT(\Delta)y_{n}$.
\medskip

Suppose that $f \in \mathcal{G}$ and $g\notin G$. Therefore, $g=\Delta$ and 
$S(f,\Delta y_{n})= y_{n}S(f,\Delta)$. We have that 
$\LT(h_{1})\mid \LT(S(f,\Delta y_{n})$, for some $h_{1}\in \mathcal{G}$ since 
$\mathcal{G}$ is a Gr\"obner basis for $\mathcal{I}$. If 
$h_{1}\neq \Delta y_{n}$, then $y_{n}$ does not divide $\LT(h_{1})$ and therefore 
$\LT(h_{1})\mid \LT(S(f,\Delta)$. On the other hand, if $h_{1} = \Delta y_{n}$, 
then $\LT(h_{1}) = y_{n}\LT(\Delta)$ and in this case 
$\LT(\Delta)\mid \LT(S(f,\Delta))$. Therefore, the division process modulo 
$\mathfrak{G}$ starts. Suppose that $S(f,\Delta) = \sum_{i\geq 1} h_{i}q_{i} + r$, 
where $h_{i}\in\mathfrak{G}$, $\LT(h_{i}q_{i})\leq \LT(S(f,\Delta))$ and 
$r$ is such that $r\neq 0$ and $\LT(r)$ is not 
divisible by any element of $\LT(\mathfrak{G})$. We have 
$S(f,\Delta)y_{n} = \sum_{i\geq 1} h_{i}q_{i}y_{n} + ry_{n}$. There exists 
$h\in\mathcal{G}$ such that $\LT(h)\mid \LT(ry_{n})$. If $h = \Delta y_{n}$, then 
$\LT(h) = \LT(\Delta)y_{n}$ and it follows that $\LT(\Delta)\mid \LT(r)$, which 
is a contradiction to the fact that $\LT(r)$ is not divisible by any element 
of $\LT(\mathfrak{G})$. If $h \neq \Delta y_{n}$, then $\LT(h)\mid \LT(r)$, which 
is again a contradiction.\qed

\begin{lemma}\label{grobcor}
Suppose that $gy_{i}\in \langle g_{1},\cdots,g_{n},\Delta\rangle$, 
then $g\in \langle g_{1},\cdots,g_{n},\Delta\rangle$.
\end{lemma}

\proof Let $i=n$. We take the monomial order as above on $R$. 
Now by the Theorem \ref{grob}, $\LT(\mathfrak{G})=
\{x_{11}x_{22}\cdots x_{k-1 k-1}x_{k+i k}y_{k},\LT(\Delta)
\mid i=1,2,\cdots n-k,k= 1,2,\cdots , n-1\}$. Let us write 
$g = \sum_{i\geq 1} h_{i}q_{i} + r$, 
where $h_{i}\in\mathfrak{G}$, $\LT(h_{i}q_{i})\leq \LT(g)$ and 
$r$ is such that $r\neq 0$ and $\LT(r)$ is not 
divisible by any element of $\LT(\mathfrak{G})$. We now consider 
$gy_{n} = \sum_{i\geq 1} h_{i}q_{i}y_{n} + ry_{n}$. A similar argument as 
above leads to a contradiction. Hence, $g\in \langle g_{1},\cdots,g_{n},\Delta\rangle$.\medskip

Let $i\neq n$. We choose the lexicographic monomial order given by 
$y_{1} > \cdots > \hat{y_{i}}> \cdots > y_{n} > y_{i} > x_{ij}$, 
for all $i,j$, such that $x_{ij}>x_{i'j'}$ if $i<i'$ or if $i=i'$ and $j< j'$. We 
follow the same argument as above to prove the statement in this case.\qed
\medskip

\begin{lemma}\label{cofactor}
$\Delta y_{i}=\sum_{j=1}^{n}A_{ji}g_{j}$, where $A_{ji}$ is the cofactor of 
$x_{ji}$ in $X$.
\end{lemma}

\proof We have 
$$\Delta y_{i} = \sum_{j=1}^{n}A_{ji}x_{ji}y_{i} 
 = \sum_{j=1}^{n}A_{ji}\left(\sum_{k=1}^{n}x_{jk}y_{k}\right) 
 - \sum_{j=1}^{n}A_{ji}\left(\sum_{k\neq i}x_{jk}y_{k}\right) 
 = \sum_{j=1}^{n}A_{ji}g_{j},$$ 
 since $\sum_{j=1}^{n}A_{ji}\left(\sum_{k\neq i}x_{jk}y_{k}\right) 
 = \sum_{k\neq i}\left(\sum_{j=1}^{n}A_{ji}x_{jk}\right)y_{k} = 0$.\qed
\medskip
 
\begin{theorem}\label{primality}
$\langle g_{1},\cdots,g_{n},\Delta\rangle$ is a prime ideal.
\end{theorem}

\proof We first prove that $[g_{1},\cdots,g_{n}]=\langle g_{1},\cdots,g_{n},\Delta\rangle$. Let $g\in [g_{1},\cdots,g_{n}]$, then $gy_{n}^{e}\in \langle g_{1},\cdots,g_{n}\rangle \subseteq \langle g_{1},\cdots,g_{n},\Delta\rangle$. If $e\neq 0$, then $gy_{n}^{e-1}\in \langle g_{1},g_{2},\cdots,g_{n},\Delta\rangle$, by Lemma \ref{grobcor}. By a repeated application of Lemma \ref{grobcor}, we finally get $g\in \langle g_{1},\cdots,g_{n},\Delta\rangle$. Hence, $[g_{1},\cdots,g_{n}]\subset\langle g_{1},\cdots,g_{n},\Delta\rangle$. Now 
$y_{i}\Delta \in \langle g_{1},\cdots,g_{n}\rangle$ by Lemma \ref{cofactor} and therefore $\langle g_{1},\cdots,g_{n},\Delta\rangle\subseteq [g_{1},\cdots,g_{n}]$. We have proved that $[g_{1},\cdots,g_{n}]=\langle g_{1},\cdots,g_{n},\Delta\rangle$ and therefore by Lemma 
\ref{compirr} the ideal $\langle g_{1},\cdots,g_{n},\Delta\rangle$ is a prime ideal.\qed
\bigskip

Let us now assume that $X$ is $n\times n$ generic or generic symmetric. 
\medskip

\begin{lemma}\label{minprime}
The minimal prime ideals containing $\mathcal{I}_{n}$ 
are $\langle y_{1},\cdots,y_{n}\rangle$ and  $\langle g_{1},\cdots, g_{n},\Delta\rangle$.
\end{lemma}

\proof Let $\mathfrak{P}$ be a minimal prime ideal containing $\mathcal{I}_{n}$.
By Lemma \ref{cofactor} we get $y_{i}\Delta\in \mathcal{I}_{n}$ and 
hence $y_{i}\Delta\in\mathfrak{P}$. If $\Delta\notin\mathfrak{P}$, then 
$y_{i}\in\mathfrak{P}$ for all $i$. Hence $\mathfrak{P}=\langle y_{1},\cdots,y_{n}\rangle$. 
If $\Delta\in\mathfrak{P}$, then $\mathfrak{P}=\langle g_{1},\cdots, g_{n},\Delta\rangle$.\qed
\medskip

The above theorem tells us that 
$\sqrt{\mathcal{I}_{n}}
=\langle y_{1},\cdots,y_{n}\rangle\cap \langle g_{1},\cdots, g_{n},\Delta\rangle$. 
We now show that $\sqrt{\mathcal{I}_{n}}=\mathcal{I}_{n}$ in the following Theorem.
\medskip

\begin{lemma}\label{radical}
$\sqrt{\mathcal{I}_{n}}=\mathcal{I}_{n}$.
\end{lemma}

\proof Let $f^{k}\in \mathcal{I}_{n}$, for some $k$. We show that $f\in \mathcal{I}_{n}$. 
We know 
that $\LT(h)\mid \LT(f^{k})$, for some $h\in\mathcal{G}$. Therefore, $\LT(h)\mid \LT(f)$, 
since the leading term of $h$ is square free. We now write $f = \sum_{i\geq i}q_{i}h_{i}+r$, 
such that $h_{i}\in\mathcal{G}$ and $r$ is the remainder modulo $\mathcal{G}$. If $r=0$ 
then $f\in\mathcal{I}_{n}$ and we are done. Suppose that $r\neq 0$. We note that 
$r = f -  \sum_{i\geq i}q_{i}h_{i}\in \sqrt{\mathcal{I}_{n}}$. The elements of the 
Gr\"{o}bner basis $\mathcal{G}$ have square-free leading terms and therefore there 
must exist some $h\in \mathcal{G}$ such that $LT(h)\mid LT(r)$, contradicting the 
fact that $r$ is the remainder. Hence $r=0$.\qed 
\medskip

We can apply Lemmas \ref{minprime} and \ref{radical} to prove the primary decomposition 
of the ideal $\mathcal{I}_{n}$, as in the part (1), statement (iii) of Theorem 
\ref{primdec1}. We can follow similar steps to prove the primary decomposition 
of the ideal $\mathcal{I}_{n+1}$, as in the part (2) of Theorem \ref{primdec1}.
\bigskip

\begin{definition}
Let $R$ be a Noetherian ring and $I$ an ideal. Then $I$ is normally torsionfree if $\mathrm{Ass} R/I=\mathrm{Ass} R/I^{n}$ for $n\geq 1$.
\end{definition}
\medskip

\begin{theorem}\label{torfree}
Let $X$ denote an $m\times n$ matrix such that its entries belong to the ideal 
$\langle \{x_{ij}; \, 1\leq i \leq m, \, 1\leq j \leq n\}\rangle$, and $Y$ is $n\times 1$ generic matrix with its entries $y_{1},\ldots,y_{n}$.  
\begin{enumerate}
\item[(i)] Let $m=n$ and $X = (x_{ij})$ be generic or generic symmetric. 
\item[(ii)] Let $m=n+1$ and $X = (x_{ij})$ be generic. 
\end{enumerate}
then $I_{1}(XY)$ is normally torsionfree.
\end{theorem}
\proof Since prime ideals in the set $\mathrm{Ass} R/I_{1}(XY)$ are minimal (by theorem \ref{primdec1}), and $I_{1}(XY)=\sqrt{(I_{1}(XY))^{n}}$ for $n\geq 1$, we have   $\mathrm{Ass} R/I_{1}(XY)=\mathrm{Ass} R/(I_{1}(XY))^{n}$ for $n\geq 1$.\qed
\medskip

\begin{proposition}\label{normal}
Let $R$ be a regular local ring and $I$ a reduced ideal. If $I$ is normally torsionfree, then $I$ is normal.
\end{proposition}

\proof See proposition $1.54$ \cite{vas}.\qed
 \medskip
 
 \begin{corollary}\label{normality}
Let $R=K[[x_{ij}, y_{j}]]$ denotes the power series ring over $K$. Let $X$ denote an $m\times n$ matrix such that its entries belong to the ideal 
$\langle \{x_{ij}; \, 1\leq i \leq m, \, 1\leq j \leq n\}\rangle$, and $Y$ is $n\times 1$ generic matrix with its entries $y_{1},\ldots,y_{n}$.  
\begin{enumerate}
\item[(i)] Let $m=n$ and $X = (x_{ij})$ be generic or generic symmetric. 
\item[(ii)] Let $m=n+1$ and $X = (x_{ij})$ be generic. 
\end{enumerate}
 then $I_{1}(XY)$ is normal.
 \end{corollary} 
  
\proof Easily follows from Theorem \ref{torfree} and Proposition \ref{normal}.\qed
\bigskip 
 
\bibliographystyle{amsalpha}

\begin{thebibliography}{A}

\bibitem{fer}
{M. Ferrero, Prime Ideals in Polynomial Rings in Several Indeterminates, 
{\em Proceedings of AMS} 125(1)(1997) 67-74.
}

\bibitem{herzog}
{J. Herzog, Certain Complexes Associated to a Sequence and a Matrix, 
{\em Manuscipta Math.} 12(1974) 217--248.
}

\bibitem{hocheagon}
{M. Hochster, J. Eagon, Cohen-Macaulay rings, Invariant Theory and 
the Generic perfection of determinantal loci, {\em Amer. J. Math.} 93(4)(1971), 
1020-1058.
}

\bibitem{hunu1}
{C. Huneke, B. Ulrich, Divisor Class groups and deformations, 
{\em Amer. J. Math} Vol. 107(6)(1985) 1265-1303.
}

\bibitem{hunu2}
{C. Huneke, B. Ulrich, Residual intersections, {\em J. reine angew. Math.} 390(1988), 1-20.
}


\bibitem{johnson}
{M.R., Johnson, J. McLoud-Mann, On equations defining Veronese Rings, 
{\em Arch. Math. (Basel)} 86(3)(2006) 205-210.
}

\bibitem{sst}
{J. Saha, I. Sengupta, G. Tripathi, Ideals of the form $I_{1}(XY)$, 
{\em arXiv:1609.02765 [math.AC]} (to appear in the Journal of Symbolic Computation (2018), https://doi.org/10.1016/j.jsc.2018.06.011).
}
\bibitem{vas}
{W.Vasconcelos, Integral Closure, 
{\em Springer Monograph in Mathematics} 2005.
}
\end{thebibliography}

\end{document}